# ON THE ORTHOGONAL POLYNOMIALS ASSOCIATED WITH A LÉVY PROCESS[1]

By Josep Lluís Solé and Frederic Utzet

*Universitat Autònoma de Barcelona*

Let $X = \{X_t, t \geq 0\}$ be a càdlàg Lévy process, centered, with moments of all orders. There are two families of orthogonal polynomials associated with $X$. On one hand, the Kailath–Segall formula gives the relationship between the iterated integrals and the variations of order $n$ of $X$, and defines a family of polynomials $P_1(x_1)$, $P_2(x_1, x_2), \ldots$ that are orthogonal with respect to the joint law of the variations of $X$. On the other hand, we can construct a sequence of orthogonal polynomials $p_n^\sigma(x)$ with respect to the measure $\sigma^2\delta_0(dx) + x^2\,\nu(dx)$, where $\sigma^2$ is the variance of the Gaussian part of $X$ and $\nu$ its Lévy measure. These polynomials are the building blocks of a kind of chaotic representation of the square functionals of the Lévy process proved by Nualart and Schoutens. The main objective of this work is to study the probabilistic properties and the relationship of the two families of polynomials. In particular, the Lévy processes such that the associated polynomials $P_n(x_1, \ldots, x_n)$ depend on a fixed number of variables are characterized. Also, we give a sequence of Lévy processes that converge in the Skorohod topology to $X$, such that all variations and iterated integrals of the sequence converge to the variations and iterated integrals of $X$.

**1. Introduction.** Let $X = \{X_t, t \geq 0\}$ be a semimartingale with $X_0 = 0$. Define the iterated integrals by the recurrence

$$(1.1)\qquad P_t^{(0)} = 1, \qquad P_t^{(1)} = X_t, \ldots, P_t^{(n)} = \int_0^t P_{s-}^{(n-1)}\, dX_s$$

Received August 2006; revised August 2006.
[1]Supported by Grant MTM2006–06427 of the Ministerio de Educación y Ciencia and FEDER.

and consider the sequence of the *variations* of $X$,

$$X_t^{(1)} = X_t, \qquad X_t^{(2)} = [X,X]_t, \qquad X_t^{(n)} = \sum_{0<s\le t} (\Delta X_s)^n, \qquad n \ge 3, \tag{1.2}$$

where $\Delta X_s = X_s - X_{s-}$. The Kailath–Segall formula (see Segall and Kailath [7] or Meyer [12]) gives the relationship between $P_t^{(n)}$ and $X_t^{(n)}$:

$$P_t^{(n)} = \frac{1}{n}(P_t^{(n-1)} X_t^{(1)} - P_t^{(n-2)} X_t^{(2)} + \cdots + (-1)^{n+1} P_t^{(0)} X_t^{(n)}). \tag{1.3}$$

We deduce that $P_t^{(n)}$ is a polynomial in $X_t^{(1)}, \ldots, X_t^{(n)}$, called the *Kailath–Segall polynomial* of order $n$. Denote this polynomial by $P_n(x_1, \ldots, x_n)$, so

$$P_t^{(n)} = P_n(X_t^{(1)}, \ldots, X_t^{(n)}).$$

The explicit expression of $P_n(x_1, \ldots, x_n)$ is

$$P_n(x_1, \ldots, x_n) = (-1)^n \sum \prod_{j=1}^n \frac{(-x_j)^{m_j}}{j^{m_j} m_j!},$$

where the summation is over all nonnegative integers $m_1, \ldots, m_n$ such that $\sum_{j=1}^n j m_j = n$ (see Avram and Taqqu [3]). The polynomials $P_n(x_1, \ldots, x_n)$, $n \ge 1$, are also a particular case of generalized Appell polynomials (see Anshelevich [1]). The first three of these polynomials are

$$\begin{aligned} P_1(x_1) &= x_1, \\ P_2(x_1, x_2) &= \tfrac{1}{2} x_1^2 - \tfrac{1}{2} x_2, \\ P_3(x_1, x_2, x_3) &= \tfrac{1}{6} x_1^3 - \tfrac{1}{2} x_1 x_2 + \tfrac{1}{3} x_3. \end{aligned}$$

If $X$ is a martingale with predictable quadratic variation $\langle X, X\rangle_t = Ct$, with finite moments of all orders, then $P^{(n)}$ and $X^{(n)}$ also have moments of all orders and the iterated integrals of different order are orthogonal, that is,

$$\mathbb{E}[P_t^{(n)} P_t^{(m)}] = \frac{1}{n!} C^n t^n \delta_{nm}, \tag{1.4}$$

where $\delta_{nm} = 1$ if $n = m$ and 0 otherwise, and $C = \mathbb{E}[X_1^2]$. The orthogonality of the Kailath–Segall polynomials with respect to the law of $(X_t, X_t^{(2)}, \ldots)$ follows. This is true, in particular, for a centered Lévy process with moments of all orders.

Consider a centered Lévy process $X$ with moments of all orders and let $\sigma^2$ be the variance of its Gaussian part and $\nu$ its Lévy measure. The measure $\gamma^\sigma(dx) = \sigma^2 \delta_0(dx) + x^2 \nu(dx)$ [we also write $\gamma(dx) = x^2 \nu(dx)$] is finite and we can construct a (finite or infinite) sequence of orthogonal polynomials $p_n^\sigma(x)$ [resp. $p_n(x)$] with respect to $\gamma^\sigma$ (resp. $\gamma$). These determine a sequence of strongly orthogonal martingales related to the Teugels martingales (see Nualart and Schoutens [13]) that are the building blocks of a kind of chaotic representation of the square functionals of the Lévy process. Therefore, we will call $\{p_n^\sigma(x), n \geq 0\}$ [or $p_n(x)$, when $\sigma = 0$] the T*eugels polynomials* associated with $X$.

The main objective of this paper is to study the probabilistic properties of, and the relationship between, the polynomials $P_n(x_1, \ldots, x_n)$ and $p_n^\sigma(x)$. When working on that problem, we found three results that we think are interesting in themselves. The first one is a proof that for a general semimartingale $X$, the Doleans exponential $\mathcal{E}(uX_t)$ (fixed $t$ and $\omega$) is analytic in a certain neighborhood of the origin and that the iterated integrals are the Taylor coefficients. This part is based on a paper of Lin [11], where the result for a Lévy process was implicitly proven, and on Yablonski [19], where a generating function of the Kailath–Segall polynomials of a Lévy process was introduced.

The second result is related to the Kailath–Segall polynomials that are expressible as polynomials of a fixed set of variables. A very interesting property of the Kailath–Segall polynomials is that when you impose some restriction on the variables $x_1, x_2, \ldots$, you get different well-known families of polynomials. For example, when $X$ is a Brownian motion, then $X_t^{(2)} = t$ and $X_t^{(n)} = 0, n \geq 3$, showing that it is enough to consider the polynomials $P_n(x, t, 0, \ldots, 0)$, and it turns out that

$$P_n(x, t, 0, \ldots, 0) = H_n(x, t),$$

where $H_n(x, t)$ are the (generalized) Hermite polynomials. In a similar way, considering a compensated Poisson process, the (generalized) Charlier polynomials $C_n(x, t)$ are obtained. It is known that the Brownian motion and the compensated Poisson process are the unique Lévy processes such that the Kailath–Segall polynomials (that means, the iterated integrals) can be written as polynomials in $x$ and $t$ (see Section 3). So, a natural question is how to characterize the Lévy processes with a similar property for a finite number of variables. The answer is that they are the Lévy processes such that the Lévy measure has finite support. This is not surprising, given the paper of Sengupta and Darkar [16], where a similar result was obtained in relation to space-time harmonic polynomials. The key to our proof is that, under the appropriate conditions, only the application of linear functions to a Lévy process gives rise to another Lévy process

Finally, the third result that we would like to mention is that, under the appropriate hypothesis, it is possible to give a sequence of simple Lévy processes $\{X_k, k \geq 1\}$ that converges in the Skorohod topology to $X$ and these processes satisfy the conditions of Avram [2] in order that all variations and iterated integrals of $X_k$ converge to the variations and iterated integral of $X$. This approximating sequence is constructed using the Gauss–Jacobi mechanical quadrature formula.

**2. Doleans exponential and Kailath–Segall polynomials.** This section is inspired by the works of Lin [11] and Yablonski [19]. First, in the paper of Lin [11], it is implicit that the iterated integrals of a Lévy process are the Taylor coefficients of the Doleans exponential at the origin, a property suggested by Meyer [12], page 318, for a semimartingale. Second, Yablonski [19] introduced a generating function in order to study a family of polynomials associated with a Lévy process that turn out to be the Kailath–Segall polynomials. However, that generating function is deterministic and Yablonski gives no probabilistic interpretation of it. Here, we combine both approaches for a general semimartingale and prove that the Yablonski generating function is the Doleans exponential of the semimartingale for fixed $\omega$, which is analytic in a neighborhood of the origin, and that the Taylor coefficients are the iterated integrals. Therefore, we prove the general claim of Meyer [12].

To begin with, for the sake of easy reference, in the next remark, we collect some results obtained by Yablonski [19].

REMARK 2.1. Given a sequence of real numbers $\overline{x} = (x_1, x_2, \ldots)$ such that $\limsup_n |x_n|^{1/n} = \lambda < \infty$, Yablonski [19] defines the generating function

$$F(u, \overline{x}) = \exp\left\{\sum_{n=1}^{\infty} \frac{(-1)^{n+1}}{n} u^n x_n\right\},$$

which is analytic for $u \in (-1/\lambda, 1/\lambda)$. Yablonski proves that in the expansion

$$F(u, \overline{x}) = \sum_{n=0}^{\infty} u^n \overline{P}_n(\overline{x}), \tag{2.1}$$

the function $\overline{P}_n(\overline{x})$ is a polynomial in $x_1, \ldots, x_n$, which satisfies

$$\begin{aligned}\overline{P}_n&(x_1, \ldots, x_n)\\ &= \frac{1}{n}\big(\overline{P}_{n-1}(x_1, \ldots, x_{n-1})x_1\\ &\qquad - \overline{P}_{n-2}(x_1, \ldots, x_{n-2})x_2 + \cdots + (-1)^{n+1}\overline{P}_0 x_n\big),\end{aligned}$$

where $\overline{P}_0 = 1$. Comparing this with (1.3), we deduce that $P_n = \overline{P}_n$.

Further, Yablonski [19] also points out the following very useful properties:

$$P_n(ax_1, \ldots, a^n x_n) = a^n P_n(x_1, \ldots, x_n)$$

and

$$P_n(x_1 + y_1, \ldots, x_n + y_n) = \sum_{k=0}^{n} P_k(x_1, \ldots, x_k) P_{n-k}(y_1, \ldots, y_{n-k}). \tag{2.2}$$

Let $X=\{X_t, t\geq 0\}$ be a semimartingale with $X_0=0$. For $u\in\mathbb{R}$, consider the Doleans equation

$$Z_t = 1 + u\int_0^t Z_{s^-}\,dX_s,$$

which has a unique solution (semimartingale) given by the Doleans exponential of $uX_t$,

$$\mathcal{E}(uX_t) = \exp\{uX_t - \tfrac{1}{2}u^2\langle X^c, X^c\rangle_t\}\prod_{0<s\leq t}(1+u\Delta X_s)e^{-u\Delta X_s}. \tag{2.3}$$

Fixing $\omega$ and $t\in\mathbb{R}_+$, it is clear that if $\{X_s, s\in[0,t]\}$ is continuous or has a finite number of jumps, then $\mathcal{E}(uX_t)$ is analytic for $u\in\mathbb{R}$. The proposition below provides a general result in this direction.

PROPOSITION 2.2. *Fix $\omega\in\Omega$ (out of a set of probability zero) and $t\in\mathbb{R}_+$. There then exists $u_0\geq 1$, depending on $\omega$ and $t$, such that the function $\mathcal{E}(uX_t)$ is analytic in $u\in(-u_0,u_0)$ and*

$$\mathcal{E}(uX_t) = \sum_{n=0}^{\infty} u^n P_t^{(n)}. \tag{2.4}$$

PROOF. Fix $\omega\in\Omega$. From the expression (2.3), it follows that we only need to prove that $\prod_{0<s\leq t}(1+u\Delta X_s)e^{-u\Delta X_s}$ is analytic. Decompose this product in the following way:

$$\begin{aligned}\prod_{0<s\leq t}&(1+u\Delta X_s)e^{-u\Delta X_s}\\ &= \underbrace{\prod_{\substack{0<s\leq t\\ |\Delta X_s|<1}}(1+u\Delta X_s)e^{-u\Delta X_s}}_{(*)}\ \underbrace{\prod_{\substack{0<s\leq t\\ |\Delta X_s|\geq 1}}(1+u\Delta X_s)e^{-u\Delta X_s}}_{(**)}.\end{aligned}$$

The term $(**)$ is analytic, since there is only a finite number of factors. On the other hand, for $u\in(-1,1)$, the expression $(*)$ is positive and, taking logarithms, we have

$$\begin{aligned}\log(*) &=_{(a)} \sum_{\substack{0<s\leq t\\ |\Delta X_s|<1}}\sum_{n=2}^{\infty}\frac{(-1)^{n+1}}{n}u^n(\Delta X_s)^n\\ &=_{(b)} \sum_{n=2}^{\infty}\frac{(-1)^{n+1}}{n}u^n\sum_{\substack{0<s\leq t\\ |\Delta X_s|<1}}(\Delta X_s)^n,\end{aligned}$$

where $(a)$ follows from $\log(1+y)-y=\sum_{n=2}^\infty \frac{(-1)^{n+1}}{n}y^n$, the series being absolutely convergent for $|y|<1$ and $(b)$ is due to Fubini's theorem, which can be applied since

$$\sum_{\substack{0<s\le t\\|\Delta X_s|<1}}\sum_{n=2}^\infty\left|\frac{(-1)^{n+1}}{n}u^n(\Delta X_s)^n\right| = \sum_{n=2}^\infty \frac{|u|^n}{n}\sum_{\substack{0<s\le t\\|\Delta X_s|<1}}|\Delta X_s|^n$$
$$\le C\sum_{\substack{0<s\le t\\|\Delta X_s|<1}}(\Delta X_s)^2\le C[X,X]_t<\infty.$$

In a similar way, it is computed that

$$\limsup_n \left|\frac{(-1)^{n+1}}{n}\sum_{\substack{0<s\le t\\|\Delta X_s|<1}}(\Delta X_s)^n\right|^{1/n}\le \limsup_n\left(\frac{1}{n}[X,X]_t\right)^{1/n}\le 1.$$

The analyticity of $(*)$ for $u\in(-1,1)$ then follows.

To compute the coefficients of the expansion of $\mathcal{E}(uX_t)$, let $u_1=(\max_{0<s\le t}|\Delta X_s|)^{-1}$. Then, for $u\in(-u_1,u_1)$, $|u\Delta X_s|<1, \forall s\in(0,t]$. Therefore, we can repeat the preceding proof to obtain

$$\prod_{0<s\le t}(1+u\Delta X_s)e^{-u\Delta X_s}$$
$$= \exp\left\{-\frac{1}{2}\sum_{0<s\le t}u^2(\Delta X_s)^2+\sum_{0<s\le t}\sum_{n=3}^\infty\frac{(-1)^{n+1}}{n}u^nX^{(n)}\right\}, \tag{2.5}$$
$$u\in(-u_1,u_1),$$

and

$$\limsup_n\left|\frac{(-1)^{n+1}}{n}\sum_{0<s\le t}(\Delta X_s)^n\right|^{1/n}\le\limsup_n\left|\frac{1}{n}\frac{1}{u_1^{n-2}}[X,X]_t\right|^{1/n}\le\frac{1}{u_1}.$$

Since the right-hand side of (2.5) is the generating function of Yablonski, from Remark 2.1 we deduce the expression (2.4) for $u\in(-u_1,u_1)$. If $u_1<1$, by the principle of analytic continuation, we deduce (2.4) for $u\in(-1,1)$. Finally, we take $u_0=\max\{1,u_1\}$. □

REMARK 2.3. From the preceding proof, we deduce that for $u\in(-1,1)$,

$$\mathcal{E}(uX_t)=\exp\Big\{uX_t-\tfrac{1}{2}u^2\langle X^c,X^c\rangle_t$$
$$+\int_{(0,t]\times\{0<|x|<1\}}(\log(1+ux)-ux)\,dJ(s,x)\Big\}$$

$$\times \prod_{\substack{0<s\le t\\ |\Delta X_s|\ge 1}} (1+u\Delta X_s)e^{-u\Delta X_s},$$

where

$$J(B)=\#\{t:(t,\Delta X_t)\in B\}, \qquad B\in\mathcal{B}((0,\infty)\times\mathbb{R}_0),$$

is the jump measure of $X$, where $\mathbb{R}_0=\mathbb{R}-\{0\}$ (see Jacod and Shiryaev [6], page 69). Moreover, fixing $\omega\in\Omega$ and $t>0$, for $u\in(-u_1,u_1)$ where $u_1=(\max_{0<s\le t}|\Delta X_s|)^{-1}$ (which depends on $\omega$ and $t$),

$$\begin{aligned}\mathcal{E}(uX_t)&=\exp\left\{\sum_{n=1}^{\infty}\frac{(-1)^{n+1}}{n}u^n X_t^{(n)}\right\}\\ &=\exp\left\{uX_t-\frac{1}{2}u^2\langle X^c,X^c\rangle_t+\int_{(0,t]\times\mathbb{R}_0}(\log(1+ux)-ux)\,dJ(s,x)\right\}.\end{aligned}$$

For any semimartingale $X$, we have that $(aX)^{(n)}=a^nX^{(n)}$. Further, given two semimartingales $X$ and $Y$ such that $[X,Y]=0$, we have

$$(X+Y)^{(n)}=X^{(n)}+Y^{(n)} \qquad \forall n\ge 1.$$

This can be proven from the expression for the product of two Doleans exponentials,

$$\mathcal{E}(uX_t)\mathcal{E}(uY_t)=\mathcal{E}(uX_t+uY_t+u^2[X,Y]_t),$$

formula (2.2) and formula (2.4). We summarize these formulas in the next proposition.

PROPOSITION 2.4. *Let $X$ and $Y$ be two semimartingales such that $X_0=Y_0=0$ and $[X,Y]=0$. For $a,b\in\mathbb{R}$, write $Z=aX+bY$. Then,*

$$\begin{aligned}&P_n(Z_t^{(1)},\dots,Z_t^{(n)})\\ &\qquad=\sum_{k=0}^{n}a^kb^{n-k}P_k(X_t^{(1)},\dots,X_t^{(k)})P_{n-k}(Y_t^{(1)},\dots,Y_t^{(n-k)}).\end{aligned}\tag{2.6}$$

**3. Kailath–Segall polynomials associated with a Lévy process.** From now on, consider that $X=\{X_t,t\ge 0\}$ is a Lévy process (meaning that $X$ has stationary and independent increments, is continuous in probability and has $X_0=0$), càdlàg, centered and with $\mathbb{E}[|X_1|^n]<\infty$ for every $n\ge 1$. Denote by $\sigma^2$ the variance of the Gaussian part of $X$ and by $\nu$ its Lévy measure. Since it is a martingale with predictable quadratic variation

$$\langle X,X\rangle_t=\left(\sigma^2+\int_{\mathbb{R}}x^2\nu(dx)\right)t,$$

it follows that the Kailath–Segall polynomials are orthogonal.

When $X$ is a Brownian motion $W=\{W_t, t\ge 0\}$,

$$X^{(1)}=W_t, \qquad X_t^{(2)}=t \quad \text{and} \quad X^{(n)}=0, \qquad n\ge 3,$$

and, therefore,

$$P_n(X_t^{(1)},\dots,X_t^{(n)})=P_n(W_t,t,0,\dots,0)=Q_n(W_t,t).$$

That is, it is enough to consider the Kailath–Segall polynomials with the variables

$$x_1=x, \qquad x_2=t \quad \text{and} \quad x_n=0, \qquad n\ge 3,$$

and then

$$P_n(x,t,0,\dots,0)=Q_n(x,t)=H_n(x,t),$$

where $H_n(x,t)$ are the Hermite polynomials defined via the generating function

$$\exp\{ux-\tfrac{1}{2}u^2t\}=\sum_{n=0}^\infty u^n H_n(x,t).$$

Note that the leading coefficient of $H_n(x,t)$ is $1/n!$, and for fixed $t>0$, $H_n(x,t)$ and $H_m(x,t)$, $n\ne m$, are orthogonal with respect to the Gaussian measure $\mathcal{N}(0,t)$. However, observe that $H_n(x,t)$ is an ordinary polynomial in $x$ and $t$ and is defined for all $x,t\in\mathbb{R}$. For $t=0$,

$$P_n(x,0,0,\dots)=H_n(x,0)=\frac{x^n}{n!}.$$

If $X$ is a compensated Poisson process of parameter $b>0$ and jumps size $a$, that is, $X_t=a(N_t-bt)$, where $N=\{N_t, t\ge 0\}$ is a standard Poisson process of intensity $b$, then

$$X_t^{(n)}=a^nN_t, \qquad n\ge 2.$$

Therefore,

$$P_n(X_t^{(1)},\dots,X_t^n)=P_n(X_t,aX_t+a^2bt,\dots,a^{n-1}X_t+a^nbt)=Q_n(X_t,t)$$

and the polynomial $Q_n(x,t)$ can be explicitly computed in the following way: write $x=a(y-bt)$, then

$$Q_n(a(y-bt),t)=a^nC_n(y,bt),$$

where $C_n(x,t)$ is the Charlier polynomial with leading coefficient $1/n!$ defined by

$$e^{-tu}(1+u)^x=\sum_{n=0}^\infty u^nC_n(x,t).$$

Again, note that $C_n(x,t)$ is defined for every $x,t\in\mathbb{R}$ and, in particular,

$$C_n(x,0)=P_n(x,x,\dots)=\frac{[x]_n}{n!},$$

where $[x]_n$ is the falling factorial, $[x]_n=x(x-1)\cdots(x-n+1), [x]_0=1$. Fixing $t>0$, the polynomials $C_n(x,t), n\geq 1$, are orthogonal with respect to the Poisson distribution of parameter $t$.

Moreover, it is known that the Brownian motion and the compensated Poisson process are the unique Lévy processes such that the Kailath–Segall polynomials can be written as polynomials in $x$ and $t$. This fact follows from Feinsilver [8], who gives a necessary condition for the iterated integral $P_t^{(n)}$ to be a polynomial in $X_t$, that condition being satisfied only by the binomial, negative binomial, Gamma, Poisson and Gaussian types (see Feinsilver [8], page 301). It is easy to check that $P_t^{(2)}$ is not a polynomial on $X_t$ for the binomial, negative binomial and Gamma process; see also Privault et al. [14] for a different proof.

Therefore, we can ask if there are Lévy processes such that the Kailath–Segall polynomials depend on a fixed finite set of variables. The answer is affirmative and the examples are very easy to find. From Proposition 2.4, we deduce, for the process

$$X_t=\sigma W_t+a(N_t-bt),$$

where $W$ is a Brownian motion and $N$ is a Poisson process of parameter $b>0$, independent of $W$, that the Kailath–Segall polynomials can be written as polynomials in $y_0,y_1$ and $t$, which are the convolutions of the polynomials $\sigma H_\bullet(y_0,t)$ and $aC_\bullet(y_1,bt)$ described above.

More generally, a jump diffusion Lévy process

$$X_t=\sigma W_t+J_t,$$

where $J$ is a centered compound Poisson process with only a finite number of jump sizes, has Kailath–Segall polynomials expressible in a fixed finite set of variables. Specifically, let

$$X_t=\sigma W_t+\sum_{j=1}^n a_j(N_j(t)-b_jt),$$

where $W$ is a Brownian motion, $N_j$ is a Poisson process of parameter $b_j$, the processes $W,N_1,\dots,N_n$ are independent and $a_1,\dots,a_n$ are different nonzero numbers. These kinds of processes will be called *simple Lévy processes* and will play a key role. For such processes, we will see that there is a family of polynomials $Q_m(x_1,\dots,x_{n+2})$ such that, for $m\geq n+2$,

$$P_m(X_t^{(1)},\dots,X_t^{(m)})=Q_m(X_t^{(1)},\dots,X_t^{(n+1)},t).$$

Moreover, for $m \geq n+2$, let $\mathcal{R}_t^m$ be the subspace of $\mathbb{R}^m$ given by the vectors $(x_1, \ldots, x_m)$ such that there exists $(y_0, \ldots, y_n) \in \mathbb{R}^{n+1}$ with

$$x_1 = \sigma y_0 + \sum_{j=1}^n a_j(y_j - b_j t),$$

$$x_2 = \sigma^2 t + \sum_{j=1}^n a_j^2 y_j,$$

$$x_k = \sum_{j=1}^n a_j^k y_j, \qquad k = 3, \ldots, m.$$

Then, the polynomial $P_m$ restricted to $\mathcal{R}_t^m$ is a (multiple) convolution of $\sigma H_\bullet(y_0, t)$ and $a_j C_\bullet(y_j, b_j t)$, $j = 1, \ldots, n$.

If $\sigma = 0$, then

$$P_m(X_t^{(1)}, \ldots, X_t^{(m)}) = Q_m(X_t^{(1)}, \ldots, X_t^{(n)}, t), \qquad m \geq n+1,$$

and in the expression for $Q_m$ restricted to a similar subspace as above, there is no Hermite polynomial part.

Remark 3.1. To summarize the situation, $P_n(x_1, \ldots, x_n), n \geq 1$, is a family of ordinary polynomials that can be evaluated on an arbitrary sequence of real numbers or random variables. However, the most interesting properties of $P_n(x_1, \ldots, x_n)$ appear when we consider a centered Lévy process with finite moments of all orders and apply $P_n$ on the sequence of the variations of $X : (X_t^{(1)}, X_t^{(2)}, \ldots)$. Then, by the Kailath–Segall formula (1.3),

$$P_t^{(n)} = P_n(X_t^{(1)}, \ldots, X_t^{(n)}),$$

and $P_n(X_t^{(1)}, \ldots, X_t^{(n)})$ and $P_m(X_t^{(1)}, \ldots, X_t^{(m)})$ are orthogonal if $n \neq m$.

From an equivalent point of view, let $\mathbb{R}^\infty = \{(x_1, x_2, \ldots), x_n \in \mathbb{R}\}$ be the set of sequences of real numbers and $\mathbb{P}_t^\infty$ the probability on $(\mathbb{R}^\infty, \mathcal{B}(\mathbb{R})^\infty)$ induced by $(X_t^{(1)}, X_t^{(2)}, \ldots)$. Writing $P_n(x_1, \ldots, x_m, \ldots) = P_n(x_1, \ldots, x_n)$, $P_n$ can be considered as a polynomial defined on $\mathbb{R}^\infty$ and the different polynomials are orthogonal:

$$\int_{\mathbb{R}^\infty} P_n P_m \, d\mathbb{P}_t^\infty = 0 \qquad \text{if } n \neq m.$$

In some cases, the probability $\mathbb{P}_t^\infty$ is concentrated in a finite-dimensional subspace of $\mathbb{R}^\infty$ and then the restriction of $P_n$ to this subspace gives rise to a new family of polynomials that depend on a finite set of variables.

So, a natural question is whether there are other examples, different from simple Lévy processes, where the Kailath–Segall polynomials can be written as polynomials in a finite number of variables. We will prove that the answer is "no."

3.1. *Polynomials of a Lévy process.* The purpose of this subsection is to study when a polynomial of a Lévy process can be a Lévy process. The result is given in the following proposition.

PROPOSITION 3.2. *Let $\{(Y(t), X_1(t), \ldots, X_d(t)), t \geq 0\}$ be a $d+1$-dimensional Lévy process with moments of all orders and $P(x_1, \ldots, x_d, t)$ a polynomial. If $Y(t) = P(X_1(t), \ldots, X_d(t), t)$, then $P(x_1, \ldots, x_d, t)$ has degree* 1.

In order to prove this proposition, we will need the following elementary property.

LEMMA 3.3. *Let $A$ be a $d \times d$ nonnegative definite matrix of* $\operatorname{rang}(A) = r \leq d$ *and let $f = (f_1(t), \ldots, f_d(t))'$ be a vector of real functions such that*

$$f' A f = 0.$$

*Then:*

1. *if $r = d$, then $f = 0$;*
2. *if $r < d$, then there are $d - r$ functions between $f_1, \ldots, f_d$, such that the other $r$ depend linearly on them.*

PROOF. This is a consequence of the fact that for any matrix $C$ conformable with $A$ such that $C'AC = 0$, we have $AC = 0$. □

PROOF OF PROPOSITION 3.2. The idea of the proof is that from the Itô formula, the decomposition of $Y_t$ as a special semimartingale is obtained, and Jacod and Shiryayev [6], Corollary II.4.19, give necessary and sufficient conditions in order that a semimartingale be a Lévy process. Then, we will prove that for a polynomial of degree $n$ of a Lévy process to be a Lévy process, it is necessary that a polynomial of order $n-1$ be a Lévy process. Hence, we can reduce to the case where the polynomial has degree 2. The proof is as follows.

First, since a polynomial in $(X_1(t), \ldots, X_d(t), t)$ can be written as a polynomial in $(X_1(t) - \mathbb{E}[X_1(t)], \ldots, X_d(t) - \mathbb{E}[X_d(t)], t)$, we can assume that the Lévy process is centered. Also, every linear combination $\sum_{j=1}^{d} \lambda_j X_j(t) + \mu t$ is a Lévy process (jointly with $Y$) and we can then eliminate such linear combinations from $P(X_1(t), \ldots, X_d(t), t)$; that is, we will assume that every monomial in $P(x_1, \ldots, x_d, t)$ has degree $\geq 2$.

1. *Degree of* $P(x_1, \ldots, x_d, t) = 2$. In this step, we will prove that if a quadratic form $P(X_1(t), \ldots, X_d(t), t)$ is a Lévy process, then $P(x_1, \ldots, x_d, t) \equiv 0$. For now, we write $x_{d+1} = t$. Then,

$$P(x_1, \ldots, x_{d+1}) = \sum_{i=1}^{d+1} b_i x_i^2 + \sum_{i<j} c_{ij} x_i x_j.$$

First, by a linear transformation, we can assume that

$$P(x_1,\dots,x_{d+1})=\sum_{i=1}^{d+1} b_i x_i^2.$$

We then have

$$Y_t=\sum_{j=1}^{d} b_j X_j^2(t)+b_{d+1}t^2.$$

Taking expectations, and recalling that the Lévy process is centered, we get

$$Ct=b_{d+1}t^2.$$

So, $b_{d+1}=0$. We will prove by induction over $d$ that also $b_1=\cdots=b_d=0$.

1.1. Let $d=1$. In order to prove that $X_t^2$ cannot be a Lévy process, write $Y_t=X_t^2$ and let $C_2=E[X_1^2]=\sigma^2+\int_{\mathbb{R}} x^2\nu(dx)$ be the variance and cumulant of order 2 of $X_1$ and $C_4=\int_{\mathbb{R}} x^4\nu(dx)$ be the cumulant of order 4 of $X_1$. On one side, if $Y$ were a Lévy process, then $\mathbb{E}[(Y_t-\mathbb{E}Y_t)^2]=Ct$, for some $C$. On the other side, from the relationship between moments and cumulants, we have

$$\mathbb{E}[Y_t^2]=\mathbb{E}[X_t^4]=C_4t+3C_2^2t^2.$$

Comparing both expressions for $\mathbb{E}[Y_t^2]$, we deduce that $X=0$.

1.2. Consider $d\ge 2$ and let $\mathbf{X}_t=(X_1(t),\dots,X_d(t))'$ be given by

$$\mathbf{X}_t=\mathbf{B}_t+\int_{(0,t]\times\mathbb{R}_0^d} x\widetilde{\mathbf{N}}(ds,d\mathbf{x}),$$

where $\mathbf{B}_t=(B_1(t),\dots,B_d(t))'$ is a $d$-dimensional Brownian motion with covariance matrix $\mathbf{A}$ and $\mathbf{N}(t,\mathbf{x})$ is the jump measure of the process, where $\mathbb{R}_0^d=\mathbb{R}^d-\{\mathbf{0}\}$ and $d\widetilde{\mathbf{N}}(t,\mathbf{x})=d\mathbf{N}(t,\mathbf{x})-dt\,d\boldsymbol{\nu}(\mathbf{x})$ is the compensated jump measure.

By Itô's formula,

$$\begin{aligned} P(\mathbf{X}_t)&=2\int_0^t\sum_{i=1}^d b_iX_i(s-)\,dB_i(s)\\ &\quad+\int_{(0,t]\times\mathbb{R}_0^d}\sum_{i=1}^d(2b_ix_iX_i(s-)+b_ix_i^2)\,d\widetilde{\mathbf{N}}(ds,d\mathbf{x}) \end{aligned}\tag{3.1}$$

$$\quad+t\left(\sum_{i=1}^d A_{ii}b_i+\int_{R^d}\left(\sum_{i=1}^d b_ix_i^2\right)\boldsymbol{\nu}(d\mathbf{x})\right).\tag{3.2}$$

The right-hand side of (3.1) is a martingale and (3.2) is of bounded variation and continuous, so the above expression is the decomposition of $P(\mathbf{X}_t)$ as a special semimartingale.

1.2.1. Assume $\mathbf{A} \neq 0$. By Jacod and Shiryavev [6], Corollary II.4.19, a necessary condition for $P(\mathbf{X}_t)$ to be a Lévy process is that the quadratic variation of the continuous martingale part should be of the form $Ct$ (the truncation function does not play any role in that condition). Then,

$$\left\langle \int_0^t \sum_{i=1}^d b_i X_i(s-)\,dB_i(s), \int_0^t \sum_{i=1}^d b_i X_i(s-)\,dB_i(s) \right\rangle = Ct \qquad \forall t \geq 0.$$

This implies that

$$\sum_{i,j} b_i b_j A_{i,j} X_i(t) X_j(t) = C \qquad \forall t \geq 0,$$

and from the fact $X_i(0) = 0$ for $i = 1, \dots, d$, we deduce that $C = 0$. We can write this expression in a vector form,

$$\mathbf{U}_t' \mathbf{A} \mathbf{U}_t = 0,$$

where $\mathbf{U}_t = (b_1 X_1(t), \dots, b_d X_d(t))'$. By Lemma 3.3, we see that if some $b_i \neq 0$, then there is a linear relationship between $X_1, \dots, X_d$.

1.2.2. If $\mathbf{A} = 0$, then

$$\int_{(0,t]\times\mathbb{R}_0^d} \sum_{i=1}^d (2 b_i x_i X_i(s-) + b_i^2 x_i^2)\,d\widetilde{\mathbf{N}}(ds, d\mathbf{x})$$

is a Lévy process. From the fact that $\int_{(0,t]\times\mathbb{R}_0^d} x_i^2\,d\mathbf{N}(ds, d\mathbf{x}) = [X_i, X_i]_t$, and as we are assuming that $(Y_t, X_t)$ is a Lévy process, we deduce that

$$\int_{(0,t]\times\mathbb{R}_0^d} \sum_{i=1}^d (b_i x_i X_i(s-))\,d\widetilde{\mathbf{N}}(ds, d\mathbf{x})$$

is a Lévy process and a martingale. Therefore,

$$\mathbb{E}\left[\left(\int_{(0,t]\times\mathbb{R}_0^d} \sum_{i=1}^d (b_i x_i X_i(s))\,d\widetilde{N}(ds, dx)\right)^2\right]$$
$$= \int_{(0,t]\times\mathbb{R}_0^d} \mathbb{E}\left[\left(\sum_{i=1}^d (b_i x_i X_i(s))\right)^2\right] ds\nu(dx),$$

but the left-hand side has the form $Ct$, and the right-hand side $C't^2$. Hence,

$$\sum_{i=1}^d b_i x_i X_i(t) = 0 \qquad \forall (x, t, \omega),\ \nu(dx) \otimes dt \otimes P \text{ a.e.}$$

From the fact that $\nu \neq 0$, it follows that there are $x_1, \dots, x_d$, not all 0, such that

$$\sum_{i=1}^d b_i x_i X_i(t) = 0 \qquad \forall (t, \omega),\ dt \otimes P \text{ a.e.}$$

and this implies a linear relationship between $X_1, \dots, X_d$.

1.2.3. From 1.2.1 and 1.2.2, we then deduce that if $P \not\equiv 0$, there is a linear relationship between $X_1, \dots, X_d$, and it follows that there is a Lévy process $(Y(t), X_1(t), \dots, X_{d-1}(t))$ and a polynomial of degree 2 such that $Y_t = P(X_1(t), \dots, X_{d-1}(t))$. Iterating the procedure, we arrive at the case $d = 1$, which is absurd (Point 1.1).

2. Degree of $P(x_1, \dots, x_d, t) = n \geq 3$. This proof is very similar to 1.2.1 and 1.2.2. As in point 1.2, we apply Itô's formula.

2.1. Assume that the covariance of the Gaussian part of $\mathbf{X}$ is not zero: $\mathbf{A} \neq 0$. The continuous martingale in $P(\mathbf{X}_t)$ is.

$$\int_0^t \sum_{i=1}^d \frac{\partial P}{\partial x_i}(\mathbf{X}_s, s)\, dB_i(t).$$

By Jacod–Shiryaev [6], its quadratic variation should be $Ct$. Write $\mathbf{V} = (\frac{\partial P}{\partial x_1}(\mathbf{X}_s, s), \dots, \frac{\partial P}{\partial x_d}(\mathbf{X}_s, s))'$. We then have

$$\mathbf{V}'\mathbf{A}\mathbf{V} = 0.$$

By Lemma 3.3, it follows that there are numbers $g_1, \dots, g_d$, not all null, such that

$$\sum_{i=1}^d g_i \frac{\partial P}{\partial x_i}(\mathbf{X}_s, s) = 0,$$

and the expression on the left-hand side is a polynomial of degree $n-1$ on $X_1(t), \dots, X_d(t), t$.

2.2. If $\mathbf{A} = 0$, then consider the bounded variation part of $P(\mathbf{X}_t, t)$,

$$V_t = \int_0^t \Biggl( \frac{\partial P}{\partial t}(\mathbf{X}_s, s) + \int_{\mathbb{R}^d} \Biggl( P(\mathbf{X}_{s-} + \mathbf{x}, s) - P(\mathbf{X}_{s-}, s) - \sum_{i=1}^d x_i \frac{\partial P}{\partial x_i}(\mathbf{X}_{s^-}, s) \Biggr) \boldsymbol{\nu}(d\mathbf{x}) \Biggr)\, ds,$$

which is continuous, thus predictable. On the other hand, the Lévy–Itô expression for a Lévy process also gives its decomposition as a special semimartingale, so, by the unicity of the decomposition, we deduce that $V_t = Ct$. Then,

$$\frac{\partial P}{\partial t}(\mathbf{X}_t, t) + \int_{\mathbb{R}^d} \Biggl( P(\mathbf{X}_{t-} + \mathbf{x}, t) - P(\mathbf{X}_{t-}, t) - \sum_{i=1}^d x_i \frac{\partial P}{\partial x_i}(\mathbf{X}_{t^-}, t) \Biggr) \boldsymbol{\nu}(d\mathbf{x}) = C.$$

This expression is also a polynomial of degree $n-1$ in $X_1(t), \dots, X_d(t), t$.

2.3. From 2.1 and 2.2, we deduce that for a polynomial of order $n$ in $X_1(t), \dots, X_d(t), t$, to be a Lévy process, it is necessary that a polynomial of order $n-1$ be a Lévy process. Iterating, we arrive at a contradiction with step 2. □

REMARK 3.4. An indication that the property expressed in Proposition 3.2 may be true for more general functions is the following. Instead of a polynomial, consider a general (sufficiently regular) function $f(x_1,\dots,x_d)$. Assume that the covariance matrix $A$ of the Gaussian part of $\mathbf{X}$ is nonsingular. The necessary condition of Jacod and Shiryavev [6], Corollary II.4.19, for $f(\mathbf{X}_t)$ to be a Lévy process becomes

$$\sum_{j=1}^{d}\left(\frac{\partial g}{\partial x_j}(\mathbf{X}_{t-})\right)^2 = C \qquad \text{a.s., for all } t\ge 0,$$

where $g$ is a function obtained from $f$ through linear changes of variable. Since the support of $\mathbf{X}_t$ is $\mathbb{R}^d$, it follows that

$$\|\nabla g(x_1,\dots,x_d)\|^2 = C \qquad \forall (x_1,\dots,x_d)\in\mathbb{R}^d.$$

This is the eikonal equation, which, in $\mathbb{R}^d$, has a unique solution given by a linear function; see Khavinson [9], Remark (ii), or Letac and Pradines [10].

3.2. *The $n$–dimensional variation process $(X^{(1)},\dots,X^{(n)})$.* We return to the general Lévy process $X$ with moments of all orders. From $\mathbb{E}[|X_1|^k]<\infty$, for all $k\ge 1$, it follows that $\mathbb{E}[|X_1^{(n)}|^k]<\infty$, for all $n,k\ge 1$ and

$$\mathbb{E}[X_t^{(1)}]=0, \qquad \mathbb{E}[X_t^{(2)}]=\left(\sigma^2+\int_{\mathbb{R}} x^2\nu(dx)\right)t$$

and

$$\mathbb{E}[X_t^{(n)}]=t\int_{\mathbb{R}} x^n\nu(dx), \qquad n\ge 3.$$

Consider the multivariate Lévy process $(X_t^{(1)},\dots,X_t^{(n)})$. Its Lévy measure $\boldsymbol{\nu}_n$ (on $\mathbb{R}^n$) is the image measure of $\nu$ by the application

$$\begin{aligned}\mathbb{R} &\longrightarrow \mathbb{R}^n,\\ x &\mapsto (x,x^2,\dots,x^n).\end{aligned}$$

By the image measure theorem, for $f:\mathbb{R}^n\to\mathbb{R}$ measurable, positive or $\boldsymbol{\nu_n}$-integrable,

$$\int_{\mathbb{R}^n} f(\mathbf{x})\boldsymbol{\nu}_n(d\mathbf{x}) = \int_{\mathbb{R}} f(x,x^2,\dots,x^n)\nu(dx).$$

The characteristic function of $(X_t^{(1)},\dots,X_t^{(n)})$ is then

$$\varphi_t(\mathbf{z})=\exp\left\{-\tfrac12 tz_1^2\sigma^2+itz_2\sigma^2+t\int_{\mathbb{R}}(e^{i\sum_{j=1}^n z_jx^j}-1-iz_1x)\nu(dx)\right\}, \tag{3.3}$$

where $\mathbf{z}=(z_1,\dots,z_n)$. Hence, the characteristic function of a linear combination

$$Z_t=\sum_{j=1}^n a_j X_t^{(j)}+a_{n+1}t,$$

such that $\mathbb{E}[Z_t]=0$, is

$$\begin{aligned}\phi_t(z)=\exp\biggl\{&-\tfrac{1}{2}ta_1^2z^2\sigma^2\\ &+t\int_{\mathbb{R}}\biggl(e^{iz\sum_{j=1}^n a_jx^j}-1-iz\sum_{j=1}^n a_jx^j\biggr)\nu(dx)\biggr\},\qquad z\in\mathbb{R}.\end{aligned}\tag{3.4}$$

Before proceeding to the main theorem of this section, we need the following lemma which will allow us to work with characteristic functions like (3.4).

Lemma 3.5. *Let $\nu$ be a Lévy measure on $(\mathbb{R},\mathcal{B}(\mathbb{R}))$ and let $f:\mathbb{R}\to\mathbb{R}$ be a continuous function such that $f(0)=0$ and $\int_{\mathbb{R}} f^2(x)\nu(dx)<\infty$. If*

$$\int_{\mathbb{R}}(e^{izf(x)}-1-izf(x))\nu(dx)=0\qquad\forall z\in\mathbb{R},$$

*then $f=0,\nu$-a.e.*

Proof. Let $\nu_f$ be the measure image of $\nu$ by $f$. From the hypothesis, it follows that $\nu_f$ is a Lévy measure and $\int_{\{|x|>1\}}|x|\nu_f(dx)<\infty$. Consider the infinitely divisible distribution $\Lambda$ that has Lévy generating triplet given by $\sigma=0$, Lévy measure $\nu_f$ and $\gamma=-\int_{\{|x|>1\}}x\nu_f(dx)$ (for this notation, see Sato [15], pages 39 and 163); its characteristic function is

$$\exp\biggl\{\int_{\mathbb{R}}(e^{izf(x)}-1-izf(x))\nu(dx)\biggr\}.$$

So, by hypothesis, $\Lambda=\delta_0$. Hence, $\nu_f=0$ and thus $f=0,\nu$-a.e. $\square$

Remark 3.6. Note that if a Lévy process has finite moment of order $k\ge 2$, then a polynomial of order $[k/2]$ without independent term satisfies the conditions on the function $f$ of the lemma.

3.3. *Kailath–Segall polynomials and finitely supported Lévy measures.*

Theorem 3.7. *There exists a number $k\ge 1$ and a family of polynomials $\{Q_n(x_1,\dots,x_k,t),n\ge k\}$, $Q_n$ of degree $n$, such that*

$$P_n(X_t^{(1)},\dots,X_t^{(n)})=Q_n(X_t^{(1)},\dots,X_t^{(k)},t),\qquad n\ge k,$$

*being this $k$ the minimum number that satisfies that condition, if and only if*

$$\mathbf{1}_{\{\sigma\neq 0\}} + \#\,\mathrm{Supp}(\nu) = k,$$

*where* $\mathrm{Supp}(\nu)$ *is the support of the Lévy measure* $\nu$.

PROOF. We first prove that the condition is sufficient.

*Case* 1. Let $\sigma = 0$ and $\mathrm{Supp}(\nu) = \{a_1, \dots, a_k\}$. Consider the polynomial of degree $k+1$,

$$R(x) = x\prod_{j=1}^{k}(x - a_j) = \sum_{j=1}^{k} c_j x^j + x^{k+1},$$

that satisfies $R(x) = 0$, $\nu$-a.e.

Denote by $L_R$ the polynomial of order 1 in $x_1, \dots, x_{k+1}$ defined by the coefficients of $R$,

$$L_R(x_1, \dots, x_{k+1}) = \sum_{j=1}^{k} c_j x_j + x_{k+1},$$

and let

$$c_{k+1} = -\mathbb{E}[L_R(X_1^{(1)}, \dots, X_1^{(k+1)})].$$

Then, the characteristic function of $L_R(X_t^{(1)}, \dots, X_t^{(k+1)}) + c_{k+1}t$ is [see (3.4)]

$$\exp\biggl\{t\int_{\mathbb{R}}(e^{izR(x)} - 1 - izR(x))\nu(dx)\biggr\},$$

which is equal to 1 because $R(x) = 0$, $\nu$-a.e. So,

$$X_t^{(k+1)} = -\sum_{j=1}^{k} c_j X_t^{(j)} - c_{k+1}t$$

and it follows that

$$P_{k+1}(X_t^{(1)}, \dots, X_t^{(k+1)}) = Q_{k+1}(X_t^{(1)}, \dots, X_t^{(k)}, t)$$

for some polynomial $Q_{k+1}$ of degree $k+1$. Now, observe that for every $n \geq k+1$, the linear system

$$\begin{aligned} a_1 g_1 + a_1^2 g_2 + \cdots + a_1^k g_k &= a_1^n \\ \vdots \qquad\qquad &\quad \vdots \\ a_k g_1 + a_k^2 g_2 + \cdots + a_k^k g_k &= a_k^n \end{aligned}$$

has a unique solution $g_1, \dots, g_k$ since the determinant of the system is a Vandermonde one. Hence,

$$\sum_{j=1}^{k} g_j x^j - x^n = 0, \qquad \nu\text{-a.e.}$$

Therefore,

$$X_t^{(n)} = \sum_{j=1}^{k} g_j X_t^{(j)} + d_n t,$$

where $d_n$ is defined is a similar way as before and thus

$$P_n(X_t^{(1)}, \dots, X_t^{(n)}) = Q_n(X_t^{(1)}, \dots, X_t^{(k)}, t) \qquad \forall n \geq k.$$

*Case* 2. When $\sigma \neq 0$ and $\mathrm{Supp}(\nu) = \{a_1, \dots, a_{k-1}\}$, consider the polynomial of order $k+1$ without independent and linear terms,

$$R(x) = x^2 \prod_{j=1}^{k-1} (x - a_j).$$

Working as in case 1, we have

$$X_t^{(k+1)} = -\sum_{j=2}^{k} c_j X_t^{(j)} - c_{k+1} t.$$

*Necessity of the condition.* Assume that

$$P_{k+1}(X_t^{(1)}, \dots, X_t^{(k+1)}) = Q_{k+1}(X_t^{(1)}, \dots, X_t^{(k)}, t).$$

In the left-hand side, by (1.3), the process $X^{(k+1)}$ appears to be simply just multiplied by $1/(k+1)$, so

$$X_t^{(k+1)} = \mathrm{Pol}(X_t^{(1)}, \dots, X_t^{(k)}, t),$$

where "Pol" means a polynomial in the specified variables. By Proposition 3.2,

$$X_t^{(k+1)} = \sum_{j=1}^{k} c_j X_t^{(j)} + c_{k+1} t$$

and there is no linear relationship between any of the $X^{(1)}, \dots, X^{(k)}$. Taking expectations,

$$c_{k+1} = -c_2\left(\sigma^2 + \int_{\mathbb{R}} x^2 \nu(dx)\right) - \sum_{j=3}^{k} \int_{\mathbb{R}} x^j \nu(dx) + \int_{\mathbb{R}} x^{k+1} \nu(dx).$$

Then, the characteristic function of $\sum_{j=1}^k c_j X_t^{(j)} - X_t^{(k+1)} + c_{k+1}t$ is 1. Specifically,

$$\exp\biggl\{-\frac{1}{2}c_1^2 t z^2\sigma^2 + t\int_{\mathbb{R}}\biggl(e^{iz(\sum_{j=1}^k c_j x^j - x^{k+1})} - 1 - iz\biggl(\sum_{j=1}^k c_j x^j - x^{k+1}\biggr)\biggr)\nu(dx)\biggr\} = 1.$$

If $\sigma > 0$, then $c_1 = 0$ and, by Lemma 3.5,

$$x^{k+1} - \sum_{j=2}^k c_j x^j = 0, \qquad \nu\text{-a.e.}$$

If $\sigma = 0$, then

$$x^{k+1} - \sum_{j=1}^k c_j x^j = 0, \qquad \nu\text{-a.e.}$$

For the second case (the first one is very similar), if the polynomial $x^{k+1} - \sum_{j=1}^k c_j x^j$ has only $r < k$ real, nonzero, distinct roots, then $\operatorname{Supp}(\nu) = \{a_1, \dots, a_r\}$ and, by the sufficiency proofs above,

$$X_t^{(r+1)} = \sum_{j=1}^r c_j X_t^{(j)} + c_{r+1}t,$$

which contradicts the assumption that there is no linear relationship between the $X^{(1)}, \dots, X^{(r+1)}$. So, $\#\operatorname{Supp}(\nu) = k$. $\square$

**4. Teugels polynomials associated with a Lévy process.** In this section, we will work under Nualart–Schoutens [13] conditions on the Lévy measure $\nu$, even though some definitions only need the condition that the Lévy process has finite moments of all orders, and that can be weakened to use only finite moments up to a convenient order. The Nualart–Schoutens [13] conditions can be expressed as the existence of $\lambda > 0$ such that

$$\int_{(-1,1)^c} e^{\lambda|x|}\nu(dx) < \infty. \tag{4.1}$$

This implies that

$$\int_{\{|x|>1\}} |x|\nu(dx) < \infty \quad \text{and} \quad \int_{\mathbb{R}} |x|^n \nu(dx) < \infty \qquad \forall n \geq 2,$$

so $X_t$ has moments of all orders, and the characteristic function of $X_t$ is analytic.

Consider the measures

$$\gamma(dx) = x^2\nu(dx) \quad \text{and,} \quad \text{if } \sigma > 0, \qquad \gamma^\sigma(dx) = \sigma^2\delta_0(dx) + x^2\nu(dx).$$

Since $\nu$ has moments of all orders $\geq 2$, it follows that $\gamma$ and $\gamma^\sigma$ are finite measures with finite moments of all orders, and the probabilities $\gamma/\gamma(\mathbb{R})$ and $\gamma^\sigma/\gamma^\sigma(\mathbb{R})$ have characteristic functions that are analytic in certain neighborhoods of the origin because, if we take $\rho = \lambda/2$, then

$$\int_{(-1,1)^c} e^{\rho|x|}x^2\nu(dx) \leq \frac{2}{\rho^2}\int_{(-1,1)^c} e^{2\rho|x|}\nu(dx) < \infty$$

and also $\int_{(-1,1)} e^{\rho|x|}x^2\nu(dx) < \infty$. This implies that the characteristic functions of $\gamma/\gamma(\mathbb{R})$ and $\gamma^\sigma/\gamma^\sigma(\mathbb{R})$ are determined by their moments (see Chow and Teicher [5], Propositions 8.4.4 and 8.4.6).

We can construct a (finite or infinite) sequence $p_n(x), n \geq 0$, of orthogonal monic polynomials with respect to $\gamma$ and another sequence of monic polynomials $p_n^\sigma(x)$ orthogonal with respect to $\gamma^\sigma$. By convention, it is always the case that $p_0(x) = p_0^\sigma(x) = 1$.

Examples.

1. Brownian motion. $X = W$. Then $\nu = 0$. The polynomials are $p_1(x) = x$, $p_n(x) = 0, n \geq 2$.
2. Standard Poisson process. $X = N_t - t$. Then $\sigma = 0$ and $\nu = \delta_1$, $p_1(x) = x - 1$, $p_n(x) = 0, n \geq 2$.
3. Simple Lévy process with 2 jump sizes, with $\sigma = 0$,

$$X_t = a_1(N_1(t) - b_1 t) + a_2(N_2(t) - b_2 t).$$

Then $\nu = b_1\delta_{a_1} + b_2\delta_{a_2}$ and

$$\begin{aligned} p_1(x) &= x - \frac{b_1 a_1^3 + b_2 a_2^3}{b_1 a_1^2 + b_2 a_2^2}, \\ p_2(x) &= (x - a_1)(x - a_2), \\ p_n(x) &= 0, \qquad n \geq 3. \end{aligned}$$

4. Gamma process. (Schoutens [17]) Let $\{G_t, t \geq 0\}$ be a Gamma process, that is, a Lévy process such that $G_t$ has distribution Gamma with mean $t$ and scale parameter equal to 1. Consider $X_t = G_t - t$. Then, $\sigma = 0$ and the Lévy measure is $\nu(dx) = \mathbf{1}_{\{x>0\}}\frac{e^{-x}}{x}\,dx$, which has infinite support. The sequence of orthogonal polynomials is infinite and they are the Laguerre polynomials $L_n^{(1)}(x)$.

As Schoutens [17] shows, it is straightforward to obtain $p_n^\sigma(x)$ from $p_n(x)$ through a family of kernels polynomials. However, in Corollary 4.3, we will see a useful relationship between $p_{n+1}^\sigma(x)$ and $p_n(x)$.

In order to compute $p_n(x)$, write

$$\Delta_n = \begin{vmatrix} m_2 & m_3 & \cdots & m_{n+2} \\ m_3 & m_4 & \cdots & m_{n+3} \\ \vdots & & & \vdots \\ m_{n+2} & m_{n+3} & \cdots & m_{2n+2} \end{vmatrix}, \tag{4.2}$$

where $m_k = \int_{\mathbb{R}} x^{k-2}\gamma(dx) = \int_{\mathbb{R}} x^k \nu(dx), k \geq 2$, and

$$D_n(x) = \begin{vmatrix} m_2 & m_3 & \cdots & m_{n+2} \\ m_3 & m_4 & \cdots & m_{n+3} \\ \vdots & & & \vdots \\ m_{n+1} & m_{n+2} & \cdots & m_{2n+1} \\ 1 & x & \cdots & x^n \end{vmatrix}. \tag{4.3}$$

There are two cases.

1. If the support of $\gamma$ is infinite, then, for every $n$, $\Delta_n \neq 0$ and

$$p_n(x) = (\Delta_{n-1})^{-1} D_n(x) \tag{4.4}$$

   defines an infinite sequence of orthogonal polynomials. This follows from well-known facts about orthogonal polynomials; see Chihara [4], pages 51 and 52, and Theorem 1.3.3.
2. If $\nu = \sum_{k=1}^n b_j \delta_{a_j}$, then there are just $n$ nonzero ($\gamma$-a.e.) orthogonal polynomials $p_0, \ldots, p_{n-1}$. The expression of the (monic) polynomial $p_n$ is very easily computed, as follows:

$$p_n(x) =_{(*)} \prod_{j=1}^n (x - a_j) = \Theta_n^{-1} \begin{vmatrix} 1 & a_1 & \cdots & a_1^n \\ \vdots & & & \vdots \\ 1 & a_n & \cdots & a_n^n \\ 1 & x & \cdots & x^n \end{vmatrix}, \tag{4.5}$$

   where

$$\Theta_n = \begin{vmatrix} 1 & a_1 & \cdots & a_1^{n-1} \\ \vdots & & & \vdots \\ 1 & a_n & \cdots & a_n^{n-1} \end{vmatrix}$$

   and $(*)$ is due to the fact that this polynomial has degree $n$ and is identically zero $\gamma$-a.e., so it is orthogonal to the first $n$ orthogonal polynomials. For $m > n$, the polynomial $p_m(x)$ is also identically zero $\gamma$-a.e. and not unique.

3. When $\sigma > 0$ and $\nu = \sum_{k=1}^n b_j \delta_{a_j}$, there are also only $n+2$ determinate polynomials $p_k^\sigma(x)$, the last one being

$$p_{n+1}^\sigma(x) = x \prod_{j=1}^n (x - a_j). \tag{4.6}$$

This expression is also deduced from the fact that this polynomial satisfies $p_{n+1}^\sigma \equiv 0, \gamma^\sigma$-a.e.

The orthogonal polynomials $p_n^\sigma(x)$, or $p_n(x)$ when $\sigma = 0$, determine a sequence of strongly orthogonal normal martingales related to the Teugels martingales (see Nualart and Schoutens [13]) and that is why we call them the Teugels polynomials associated with $X$. In Section 4.3, we provide an explicit expression for those martingales.

Remark 4.1. We have changed the notation of Nualart and Schoutens [13] and Schoutens [17] because they write $p_n(x)$ for the orthogonal polynomial of degree $n-1$, and here it denotes the polynomial of degree $n$.

The next theorem is a modification of the Gauss–Jacobi mechanical quadrature formula (see Szegö [18], Theorems 3.4.1 and 3.4.2).

Theorem 4.2. *Assume that the Lévy measure $\nu$ has infinite support and let $n \geq 1$ be such that $p_n(0) \neq 0$. There are then different nonzero numbers $a_1, \ldots, a_n$ and strictly positive numbers $b_1, \ldots, b_n$ such that the (Lévy) measure with finite support,*

$$\nu_n = \sum_{k=1}^n b_k \delta_{a_k},$$

*satisfies*

$$\int_{\mathbb{R}} x^k \nu(dx) = \int_{\mathbb{R}} x^k \nu_n(dx), \qquad k = 2, \ldots, 2n+1. \tag{4.7}$$

*Moreover, let $\gamma_n(dx) = x^2 \nu_n(dx)$ and $\gamma_n^\sigma(dx) = \sigma^2 \delta_0(dx) + x^2 \nu_n(dx)$. Then, $\gamma$ and $\gamma_n$ (resp. $\gamma^\sigma$ and $\gamma_n^\sigma$) have the same orthogonal polynomials up to order $n$.*

By the Gauss–Jacobi formula (Szegö [18], Theorem 3.4.1), the numbers $a_1, \ldots, a_n$ are the $n$ different nonzero real roots of $p_n(x)$, and $b_1, \ldots, b_n$ are the unique solution of the compatible system

$$\begin{array}{ccccccccc} m_2 & = & a_1^2 b_1 & + & a_2^2 b_2 & + & \cdots & + & a_n^2 b_n \\ \vdots & & & & & & & & \vdots \\ m_{2n+1} & = & a_1^{2n+1} b_1 & + & a_2^{2n+1} b_2 & + & \cdots & + & a_n^{2n+1} b_n, \end{array}$$

where $m_k = \int_{\mathbb{R}} x^k \nu(dx), k \geq 2$. The numbers $b_1, \dots, b_n$, called the *Christoffel numbers*, are all strictly positive (Szegö [18], Theorem 3.4.2). From (4.7), it follows that the finite measures $\gamma$ and $\gamma_n$ have the same *moments* of order $0, 1, \dots, 2n-1$, denoted by $m_2, \dots, m_{2n+1}$. Hence, from the expressions (4.2), (4.3) and (4.4), we deduce that they have the same Teugels polynomials up to order $n$.

Corollary 4.3. *Let $\nu$ be a Lévy measure such that $p_n(0) \neq 0$. There are then $\lambda_n$ and $\lambda_{n-1}$ such that $p^\sigma_{n+1}(x) = xp_n(x) - \lambda_n p^\sigma_n(x) - \lambda_{n-1} p^\sigma_{n-1}(x)$.*

Proof. The polynomials $p^\sigma_{n+1}(x)$ and $xp_n(x)$ are monic and have degree $n+1$. The polynomial $xp_n(x) - p^\sigma_{n+1}(x)$ can then be written as

$$xp_n(x) - p^\sigma_{n+1}(x) = \sum_{j=0}^{n} \lambda_j p^\sigma_j(x),$$

where

$$\begin{aligned}\lambda_j &= K_j^{-1} \int_{\mathbb{R}} (xp_n(x) - p^\sigma_{n+1}(x)) p^\sigma_j(x) \gamma^\sigma(dx) \\ &= K_j^{-1} \int_{\mathbb{R}} xp_n(x) p^\sigma_j(x) \gamma^\sigma(dx)\end{aligned}$$

and

$$K_j = \int_{\mathbb{R}} (p^\sigma_j(x))^2 \gamma^\sigma(dx).$$

Consider the discrete measures $\gamma_n$ and $\gamma^\sigma_n$ of the preceding theorem. Then:

1. the measures $\gamma^\sigma$ and $\gamma^\sigma_n$ have the same moments up to order $2n-1$ and, for $j = 0, \dots, n-2$,

$$\lambda_j = K_j^{-1} \int_R xp_n(x) p^\sigma_j(x) \gamma^\sigma_n(dx).$$

2. denoting by $\widetilde{p}_j(x)$ [resp. $\widetilde{p}^\sigma_j(x)$] the Teugels polynomials of $\gamma_n$ (resp. $\gamma^\sigma_n$), we have

$$\widetilde{p}_j(x) = p_j(x) \quad \text{and} \quad \widetilde{p}^\sigma_j(x) = p^\sigma_j(x), \qquad j = 1, \dots, n,$$

and [(4.6)]

$$\widetilde{p}^\sigma_{n+1}(x) = x\widetilde{p}_n(x) = xp_n(x) = 0, \qquad \gamma^\sigma_n\text{-a.e.},$$

so it follows that, for $j = 0, \dots, n-2$,

$$\begin{aligned}\lambda_j &= K_j^{-1} \int_R xp_n(x) p^\sigma_j(x) \gamma^\sigma_n(dx) \\ &= K_j^{-1} \int_R \widetilde{p}^\sigma_{n+1}(x) \widetilde{p}^\sigma_j(x) \gamma^\sigma_n(dx) = 0.\end{aligned}$$

□

4.1. *An approximating sequence of simple Lévy processes.* An interesting consequence of Theorem 4.2 is that it provides a way to construct a sequence of simple Lévy processes that converges in the Skorohod topology to $X$, satisfying the conditions of Avram [2] in order that all variations and iterated integrals of the sequence converge to the variations and iterated integral of the limit. From the separation of zeros theorem of $p_n$ (see [4]), if $p_n(0) = 0$, then $p_{n+1}(0) \neq 0$. There is then a sequence $m_1 < m_2 < \cdots \nearrow \infty$ such that $p_{m_k}(0) \neq 0, \forall k$. Let $X_k = \{X_k(t), t \in \mathbb{R}\}$ be a centered Levy process with diffusion coefficient $\sigma$ and Lévy measure $\nu_{m_k}$ given in Theorem 4.2. That is, the law of $X_k$ is

$$X_k(t) =_{\text{Law}} \sigma W_t + \sum_{j=1}^{m_k} a_j (N_j(t) - b_j t).$$

Denote by $P_k^{(n)}$ and $X_k^{(n)}$ the iterated integral and the variation of order $n$ of $X_k$, respectively, and by $P^{(n)}$ and $X^{(n)}$ the iterated integral and variation of $X$, respectively.

THEOREM 4.4. *Let $X$ be a Lévy process that satisfies the condition (4.1) and such that $\nu$ has infinite support. With the above notation, for every $n$,*

$$\lim_k P_k^{(n)} = P^{(n)} \quad \textit{and} \quad \lim_k X_k^{(n)} = X^{(n)}$$

*(both convergences in the Skorohod sense).*

PROOF. By Avram [2], it suffices to prove that

$$\lim_k (X_k, [X_k, X_k]) = (X, [X, X]) \qquad \text{in the Skorohod sense.}$$

Since all of the process involved are Lévy process, by Jacod and Shiryaev [6], Corollary VII.3.6, it is sufficient to prove that

$$\lim_k (X_k(1), [X_k, X_k]_1) = (X(1), [X, X]_1) \qquad \text{in distribution,}$$

and by the Cramer–Wold device, this is equivalent to proving that for every $u, v \in \mathbb{R}$,

$$\lim_k (uX_k(1) + v[X_k, X_k]_1) = uX(1) + v[X, X]_1 \qquad \text{in distribution.}$$

From (3.4), the characteristic function of $uX(1) + v[X, X]_1$ is

$$\psi(z) = \exp\biggl\{-\tfrac{1}{2}u^2 z^2 \sigma^2 + izv\biggl(\sigma^2 + \int_{\mathbb{R}} x^2 \nu(dx)\biggr) + t\int_{\mathbb{R}} (e^{iz(ux+vx^2)} - 1 - iz(ux + vx^2))\nu(dx)\biggr\}.$$

From the fact that the characteristic function of $X_1$ is analytic, it follows that $\psi(z)$ also is. So, it suffices to show that all cumulants of $uX_k(1) + v[X_k, X_k]_1$ converge to the corresponding cumulants of $uX(1) + v[X, X]_1$ and this is clear from the construction of $\nu_{m_k}$. $\square$

4.2. *The relationship between Kailath–Segall polynomials and Teugels polynomials.*

4.2.1. *Preliminary results.* This subsection is purely algebraic; later, we will give a probabilistic interpretation of the results. First, it is convenient to introduce a new notation. Given a polynomial of order $n$,

$$P(x) = c_0 + c_1 x + \cdots + c_n x^n,$$

we denote by $\mathcal{L}(P)(x_1, \dots, x_{n+1})$ the polynomial of degree 1 in $x_1, \dots, x_{n+1}$ associated with the coefficients of $P$:

$$\mathcal{L}(P)(x_1, \dots, x_{n+1}) = c_0 x_1 + \cdots + c_n x_{n+1}. \tag{4.8}$$

Of course, we can recover $P(x)$ from $\mathcal{L}(P)(x_1, \dots, x_{n+1})$:

$$P(x) = \mathcal{L}(P)(1, x, \dots, x^n).$$

Second, we need to consider some finite-dimensional vector spaces. Let $\mathbf{a} = (a_1, \dots, a_n)$, where $a_1, \dots, a_n$ are different nonzero numbers. Write

$$\mathcal{S}^{\mathbf{a}}_{n+1} = \Bigg\{ (x_1, \dots, x_{n+1}) \in \mathbb{R}^{n+1} : x_1 = \sum_{j=1}^n a_j y_j, \\ x_2 = \sum_{j=1}^n a_j^2 y_j, \dots, x_{n+1} = \sum_{j=1}^n a_j^{n+1} y_j, \\ \text{for some } (y_1, \dots, y_n) \in \mathbb{R}^n \Bigg\}.$$

$\mathcal{S}^{\mathbf{a}}_{n+1}$ is subspace of dimension $n$ of $\mathbb{R}^{n+1}$, and there is the projection

$$\mathbb{R}^{n+1} \longrightarrow \mathcal{S}^{\mathbf{a}}_{n+1}, \\ (x_1, \dots, x_{n+1}) \to (x_1, \dots, x_n, u_{n+1}),$$

where $u_{n+1}$ is computed as follows. By the Vandermonde determinant property, we can find $(y_1, \dots, y_n) \in \mathbb{R}^n$ such that

$$x_1 = \sum_{j=1}^n a_j y_j, \qquad x_2 = \sum_{j=1}^n a_j^2 y_j, \dots, x_n = \sum_{j=1}^n a_j^n y_j. \tag{4.9}$$

We then write

$$u_{n+1} = \sum_{j=1}^{n} a_j^{n+1} y_j.$$

LEMMA 4.5. *With the above notation,*

$$u_{n+1} = -\mathcal{L}(P)(x_1, \ldots, x_n, 0),$$

*where*

$$P(x) = \prod_{j=1}^{n} (x - a_j).$$

PROOF. From the expression of $P(x)$ given in (4.5),

$$\mathcal{L}(P)(x_1, \ldots, x_{n+1}) = \Theta_n^{-1} \begin{vmatrix} 1 & a_1 & \cdots & a_1^n \\ \vdots & & & \vdots \\ 1 & a_n & \cdots & a_n^n \\ x_1 & x_2 & \cdots & x_{n+1} \end{vmatrix}.$$

Hence,

$$\mathcal{L}(P)(x_1, \ldots, x_n, u_{n+1}) = \mathcal{L}(P)\left(\sum_{j=1}^{n} a_j y_j, \ldots, \sum_{j=1}^{n} a_j^{n+1} y_j\right) = 0.$$

The polynomial $P(x)$ is monic, so

$$u_{n+1} = -\mathcal{L}(P)(x_1, \ldots, x_n, 0). \qquad \square$$

Define the polynomial of degree $n+1$,

$$J_{n+1}^{\mathbf{a}}(x_1, \ldots, x_n) = P_{n+1}(x_1, \ldots, x_n, u_{n+1}) = \sum \prod_{j=1}^{n} \frac{1}{h_j!} a_j^{h_j} [y_j]_{h_j},$$

where the summation is over all nonnegative integers $h_1, \ldots, h_n$ such that $\sum_{j=1}^{n} h_j = n+1$, $[x]_n$ is the falling factorial and $y_1, \ldots, y_n$ are given in (4.9). This strange expression is the multiple convolution of Charlier polynomials $aC_{\bullet}(y_j, 0), j = 1, \ldots, n$. Note that when working with the polynomials, the variable $t$ does not play the role of time and can be used freely according to our needs.

PROPOSITION 4.6. *For every* $(x_1, \ldots, x_{n+1}) \in \mathbb{R}^{n+1}$,

$$P_{n+1}(x_1, \ldots, x_{n+1}) - J_{n+1}^{\mathbf{a}}(x_1, \ldots, x_n) = \frac{(-1)^n}{n+1} \mathcal{L}(P)(x_1, \ldots, x_{n+1}),$$

*where*

$$P(x) = \prod_{j=1}^{n} (x - a_j).$$

*Equivalently,*

$$P_{n+1}(1, x, \ldots, x^n) - J_{n+1}^{\mathbf{a}}(1, x, \ldots, x^{n-1}) = \frac{(-1)^n}{n+1} P(x).$$

PROOF. Simply note that $P_{n+1}$ is linear in $x_{n+1}$, with coefficient $(-1)^n/(n+1)$ [see (1.3)], and apply Lemma 4.5. □

Note that this proposition is true if we replace $P_{n+1}$ by another polynomial linear in the variable $x_{n+1}$, but we will see that with $P_{n+1}$, it has an interesting probabilistic interpretation.

In the same way, take $\sigma > 0$ and write

$$\mathcal{S}_{n+2}^{\sigma,\mathbf{a}} = \Bigg\{ (x_1, \ldots, x_{n+2}) \in \mathbb{R}^{n+2} : x_1 = \sigma y_0 + \sum_{j=1}^{n} a_j y_j,$$
$$x_2 = \sum_{j=1}^{n} a_j^2 y_j, \ldots, x_{n+2} = \sum_{j=1}^{n} a_j^{n+1} y_j,$$
$$\text{for some } (y_0, y_1, \ldots, y_n) \in \mathbb{R}^{n+1} \Bigg\}.$$

Consider the projection

$$\mathbb{R}^{n+2} \longrightarrow \mathcal{S}_{n+2}^{\sigma,\mathbf{a}},$$
$$(x_1, \ldots, x_{n+2}) \;\to\; (x_1, \ldots, x_{n+1}, u_{n+2}^{\sigma}),$$

where

$$u_{n+2}^{\sigma} = \sum_{j=1}^{n} a_j^{n+2} y_j,$$

$y_1, \ldots, y_n$ being the solution of $x_2 = \sum_{j=1}^{n} a_j^2 y_j, \ldots, x_{n+1} = \sum_{j=1}^{n} a_j^{n+1} y_j$. With the same proof as Lemma 4.5, we have that $\mathcal{L}(P)(x_2, \ldots, x_{n+1}, u_{n+2}^{\sigma}) = 0$, where $P(x) = \prod_{j=1}^{n}(x - a_j)$. Also, note that

$$\mathcal{L}(xP)(x_1, \ldots, x_{n+2}) = \mathcal{L}(P)(x_2, \ldots, x_{n+2}).$$

Define the polynomial

$$J_{n+2}^{\sigma,\mathbf{a}}(x_1, \ldots, x_{n+1}) = P_{n+2}(x_1, \ldots, x_{n+1}, u_{n+2}^{\sigma}),$$

which has an expression similar to $J_{n+1}^{\mathbf{a}}$ with the addition of a Hermite polynomial $H_{\bullet}(y_0, 0)$. We then have the following.

PROPOSITION 4.7. *For every* $(x_1,\dots,x_{n+2})\in\mathbb{R}^{n+2}$,

$$P_{n+2}(x_1,\dots,x_{n+2}) - J^{\sigma,\mathbf{a}}_{n+2}(x_1,\dots,x_{n+1}) = \frac{(-1)^{n+1}}{n+2}\mathcal{L}(xP)(x_1,\dots,x_{n+2}),$$

*where*

$$P(x)=\prod_{j=1}^n (x-a_j).$$

*Equivalently,*

$$P_{n+2}(1,x,\dots,x^{n+1}) - J^{\sigma,\mathbf{a}}_{n+2}(1,x,\dots,x^{n}) = \frac{(-1)^{n+1}}{n+2} xP(x).$$

4.2.2. *Teugels polynomials.* The propositions of the previous subsection can be transferred when we have a Lévy measure $\nu$ and the corresponding Teugels polynomials $p_n(x)$ and $p_n^\sigma(x)$. We use Corollary 4.3 to identify these polynomials.

COROLLARY 4.8. *Fix* $n\ge 1$ *such that* $p_n(0)\neq 0$ *and let* $a_1,\dots,a_n$ *be the roots of* $p_n(x)$. *Then,*

$$P_{n+1}(x_1,\dots,x_{n+1}) - J^{\mathbf{a}}_{n+1}(x_1,\dots,x_n) = \frac{(-1)^n}{n+1}\mathcal{L}(p_n)(x_1,\dots,x_{n+1})$$

*and*

$$\begin{aligned}
&P_{n+2}(x_1,\dots,x_{n+2}) - J^{\sigma,\mathbf{a}}_{n+2}(x_1,\dots,x_{n+1})\\
&\quad= \frac{(-1)^{n+1}}{n+2}(\mathcal{L}(p^\sigma_{n+1})(x_1,\dots,x_{n+2})\\
&\qquad\qquad + \lambda_n\mathcal{L}(p^\sigma_n)(x_1,\dots,x_{n+1}) + \lambda_{n-1}\mathcal{L}(p^\sigma_{n-1})(x_1,\dots,x_n)).
\end{aligned}$$

*Equivalently,*

$$P_{n+1}(1,x,\dots,x^n) - J^{\mathbf{a}}_{n+1}(1,x,\dots,x^{n-1}) = \frac{(-1)^n}{n+1}p_n(x)$$

*and*

$$\begin{aligned}
&P_{n+2}(1,x,\dots,x^{n+1}) - J^{\sigma,\mathbf{a}}_{n+2}(1,x,\dots,x^{n})\\
&\quad= \frac{(-1)^{n+1}}{n+2}(p^\sigma_{n+1}(x) + \lambda_n p^\sigma_n(x) + \lambda_{n-1}p^\sigma_{n-1}(x)).
\end{aligned}$$

4.3. *Main result.* The Teugels martingales $\{Y^{(n)}, n \geq 1\}$ (see Nualart and Schoutens [13]) are defined by $Y_t^{(n)} = X_t^{(n)} - \mathbb{E}[X_t^{(n)}], n \geq 1$. Specifically,

$$Y_t^{(1)} = X_t, \qquad Y_t^{(2)} = X_t^{(2)} - t\biggl(\sigma^2 + \int_{\mathbb{R}} x^2 \nu(dx)\biggr)$$

and

$$Y_t^{(n)} = X_t^{(n)} - t\int_{\mathbb{R}} x^n \nu(dx), \qquad n \geq 3.$$

(Nualart and Schoutens [13] write $X_t^{(2)} = \sum_{s \leq t} (\Delta X_s)^2$, instead of $[X, X]$, as we have done; however, both definitions give the same $Y^{(2)}$.) By an orthogonalization procedure, they obtain a family $\{H^{(n)}, n \geq 1\}$ of normal martingales, pairwise strongly orthogonal, that, under the hypothesis (4.1), generate all of $L^2(\Omega)$ by sums of iterated integrals. In order to strongly orthogonalize $\{Y^{(n)}, n \geq 1\}$, if $\sigma > 0$, they show that you can look for the sequence of orthogonal polynomials $p_n^\sigma(x)$ and take

$$H_t^{(n+1)} = \mathcal{L}(p_n^\sigma)(Y_t^{(1)}, \dots, Y_t^{(n+1)}),$$

and the same expression with $p_n$ replacing $p_n^\sigma$ if $\sigma = 0$.

THEOREM 4.9. *Let $X$ be a centered Lévy process with moments of all orders and fix $n \geq 1$ such that the Teugels polynomial of order $n$, $p_n(x)$, does not have a zero root. Let $a_1, \dots, a_n$ be the roots of $p_n(x)$. If $\sigma = 0$, then*

$$\begin{aligned} &P_{n+1}(X_t^{(1)}, \dots, X_t^{(n+1)}) \\ &\quad = \frac{(-1)^n}{n+1} H_t^{(n+1)} + J_{n+1}^{\mathbf{a}}(X_t^{(1)}, \dots, X_t^{(n)}) + \frac{(-1)^n}{n+1} C_n t, \end{aligned} \tag{4.10}$$

*where*

$$C_n = \int_{\mathbb{R}} x(p_n(x) - p_n(0))\nu(dx)$$

*and $J_{n+1}^{\mathbf{a}}(X_t^{(1)}, \dots, X_t^{(n)}) + \frac{(-1)^n}{n+1} C_n t$ is orthogonal to $P_1(X_t^{(1)}), \dots, P_n(X_1^{(1)}, \dots, X_t^{(n)})$.*

*If $\sigma > 0$, then*

$$\begin{aligned} P_{n+2}(X_t^{(1)}, \dots, X_t^{(n+2)}) = {}& \frac{(-1)^{n+1}}{n+2} (H_t^{(n+2)} + \lambda_n H_t^{(n+1)} + \lambda_{n-1} H_t^{(n)}) \\ & + J_{n+2}^{\sigma, \mathbf{a}}(X_t^{(1)}, \dots, X_t^{(n+1)}) + \frac{(-1)^{n+1}}{n+2} D_{n+1} t, \end{aligned}$$

*where $\lambda_n$ and $\lambda_{n-1}$ are given in Corollary 4.3 and*

$$D_{n+1} = \sigma^2 p_n(0) + \int_{\mathbb{R}} x^2 p_n(x)\nu(dx).$$

*Moreover,* $J^{\sigma,\mathbf{a}}_{n+2}(X^{(1)}_t,\dots,X^{(n+1)}_t)+\frac{(-1)^{n+1}}{n+2}D_{n+1}t$ *is orthogonal to* $P_1(X^{(1)}_t)$, $\dots,P_{n-1}(X^{(1)}_1,\dots,X^{(n-1)}_t)$.

PROOF. Since the proof for $\sigma>0$ is very similar to the case $\sigma=0$, we consider only the latter one. Formula 4.10 follows from Corollary 4.8 and

$$\begin{aligned}H^{(n+1)}_t&=\mathcal{L}(p_n)(Y^{(1)}_t,\dots,Y^{(n+1)}_t)\\&=\mathcal{L}(p_n)(X^{(1)}_t,\dots,X^{(n+1)}_t)-t\int_{\mathbb{R}}x(p_n(x)-p_n(0))\nu(dx).\end{aligned}$$

To prove the orthogonality between $J^{\mathbf{a}}_{n+1}(X^{(1)}_t,\dots,X^{(n)}_t)+\frac{(-1)^n}{n+1}C_nt$ and $P_j(X^{(1)}_1,\dots,X^{(j)}_t)$ for $j=1,\dots,n$, observe that, by definition of $J^{\mathbf{a}}_{n+1}$ and Lemma 4.5,

$$J^{\mathbf{a}}_{n+1}(X^{(1)}_t,\dots,X^{(n)}_t)+\frac{(-1)^n}{n+1}C_nt=P_{n+1}(X^{(1)}_1,\dots,X^{(n)}_t,V^{(n+1)}_t),$$

where

$$V^{(n+1)}_t=-\mathcal{L}(p_n)(X^{(1)}_t,\dots,X^{(n)}_t,0)+C_nt.$$

The idea of the proof is to construct a simple Lévy process $Z_t$ such that, for $r\le n$,

$$\begin{aligned}&\mathbb{E}[P_r(X^{(1)}_1,\dots,X^{(r)}_t)P_{n+1}(X^{(1)}_t,\dots,X^{(n)}_t,V^{(n+1)}_t)]\\&\qquad=\mathbb{E}[P_r(Z^{(1)}_1,\dots,Z^{(r)}_t)P_{n+1}(Z^{(1)}_1,\dots,Z^{(n+1)}_t)]\end{aligned}\tag{4.11}$$

and by the orthogonality of the iterated integrals of different order, the expectation on the right is zero.

With this objective, consider the Lévy process $(X^{(1)}_t,\dots,X^{(n)}_t,V^{(n+1)}_t)$ that has characteristic function [see (3.3)]

$$\varphi_X(\mathbf{z})=\exp\biggl\{t\int_{\mathbb{R}}(e^{i(\sum_{j=1}^{n+1}z_jx^j-z_{n+1}xp_n(x))}-1-ix(z_1-z_{n+1}p_n(x)))\nu(dx)\biggr\},$$

where $\mathbf{z}=(z_1,\dots,z_{n+1})$.

On the other hand, let $b_1,\dots,b_n$ the Christoffel numbers corresponding to $\nu$ given in Theorem 4.2. Define (on the same probability space or another)

$$Z_t=\sum_{j=1}^n a_j(N_n(t)-b_jt),$$

where $N_1,\dots,N_n$ are independent Poisson processes with respective intensities $b_1,\dots,b_n$. The characteristic function of $(Z^{(1)}_t,\dots,Z^{(n+1)}_t)$ is [see (3.3)]

$$\varphi_Z(\mathbf{z})=\exp\biggl\{t\int_{\mathbb{R}}(e^{i(\sum_{j=1}^{n+1}z_jx^j)}-1-ixz_1)\nu_n(dx)\biggr\},$$

where $\nu_n = \sum_{j=1}^{n} b_j \delta_{a_j}$. By Theorem 4.2, $\nu$ and $\nu_n$ have the same moments up to order $2n+1$, both $\nu$ and $\nu_n$ have the same first $n$ Teugels polynomials and $p_n(x) \equiv 0, \nu_n$-a.e. Let $j_1, \dots, j_{n+1}$ be nonnegative integers such that $\sum_{k=1}^{n+1} k j_k \leq 2n+1$. Therefore, the (joint) cumulant of order $j_1$ in the first component, order $j_2$ in the second component and so on, of the vectors $(X_t^{(1)}, \dots, X_t^{(n)}, V_t^{(n+1)})$ and $(Z_t^{(1)}, \dots, Z_t^{(n+1)})$ are the same. So, a polynomial up to degree $2n+1$ of $(X_t^{(1)}, \dots, X_t^{(n)}, V_t^{(n+1)})$ and the same polynomial of $(Z_t^{(1)}, \dots, Z_t^{(n+1)})$ have the same expectation. In particular, for $r \leq n$, we have the identity (4.11). $\square$

Example. A very simple example will help to interpret Theorem 4.9. Consider a Lévy process $X$ with $\sigma = 0$ and Lévy measure $\nu$, and let $a$ be the root of its Teugels polynomial of order 1, $p_1(x)$. Assume $a \neq 0$ and let $b$ be the solution of

$$\int_{\mathbb{R}} x^2 \nu(dx) = ba^2.$$

Let $Z_t = a(N_t - bt)$, where $N_t$ is a Poisson process of intensity $b$. $Z_t$ is then a simple Lévy process that has Lévy measure $\nu_1 = b\delta_a$. By Gauss–Jacobi Theorem 4.2, $\nu$ and $\nu_1$ have the same moments of order 2 and 3. So, $X_t$ and $Z_t$ have the same cumulants of order 2 and 3, and, since both are centered, they have the same moments of those orders. Then, on one hand, $J_2^a(X_t^{(1)}) = P_2(X_t, aX_t)$, so

$$J_2^a(X_t) - \tfrac{1}{2} t \int_{\mathbb{R}} x^2 \nu(dx) = P_2\Big(X_t, aX_t + t\int_{\mathbb{R}} x^2 \nu(dx)\Big).$$

On the other hand,

$$P_2(Z_t^{(1)}, Z_t^{(2)}) = P_2(Z_t, aZ_t + ba^2 t).$$

We then have

$$\mathbb{E}\Big[P_1(X_t)\Big(J_2^a(X_t) - \tfrac{1}{2} t \int_{\mathbb{R}} x^2 \nu(dx)\Big)\Big] = \mathbb{E}[P_1(Z_t^{(1)}) P_2(Z_t^{(1)}, Z_t^{(2)})]$$

because $P_1(X_t)(J_2^a(X_t) - \frac{1}{2} t \int_{\mathbb{R}} x^2 \nu(dx))$ is a product of a polynomial of degree 1 and a polynomial of order 2, in $X_t$, which is centered, so the expectation of that product depends only on the moments of order 2 and 3. So, Theorem 4.9 says that we have a decomposition

$$P_2(X_t^{(1)}, X_t^{(2)}) = J_2^a(X_t) - \tfrac{1}{2} t \int_{\mathbb{R}} x^2 \nu(dx) - \tfrac{1}{2} H_t^{(2)},$$

such that:

1. $J_2^a(X_t) - \frac{1}{2} t \int_{\mathbb{R}} x^2 \nu(dx)$ is orthogonal to $P_1(X_t^{(1)})$;

2. $J_2^a(X_t) - \frac{1}{2}t\int_{\mathbb{R}} x^2\nu(dx) = P_2(X_t^{(1)}, V_t^{(2)})$, where $(X_t^{(1)}, V_t^{(2)})$ is a Lévy process that has the same moments (up to order 3) as the variations $(Z_t^{(1)}, Z_t^{(2)})$ of the simple Lévy process $Z_t$.

**Acknowledgments.** The authors wish to express their thanks to José Antonio Carrillo and Xavier Mora for their fruitful comments about the eikonal equation.

Departament de Matemàtiques
Facultat de Ciències
Universitat Autònoma de Barcelona
08193 Bellaterra (Barcelona)
Spain
E-mail: jllsole@mat.uab.cat
utzet@mat.uab.cat